\documentclass[runningheads]{svjour2}
\journalname{Probability Theory and Related Fields}
\smartqed  

\usepackage{graphicx}
\usepackage{mathptmx}      
\usepackage{amsmath,amsfonts,amssymb,stmaryrd}
\usepackage[latin1]{inputenc}

\newtheorem{thm}{Theorem}
\newtheorem{lem}{Lemma}

\newtheorem{prp}{Proposition}

\newcounter{const}
\renewcommand{\theconst}{c_{\arabic{const}}}
\newcommand{\const}[1][\relax]{%
   \bgroup\refstepcounter{const}%
   \theconst\ifx#1\relax\else\label{#1}\fi\egroup}

\newcommand{\ZZ}{\mathbb Z}
\newcommand{\EE}{\mathbb E}

\newcommand{\LL}{\mathcal L}
\newcommand{\RR}{\mathcal R}
\newcommand{\NN}{\mathcal N}
\newcommand{\CC}{\mathcal C}

\newcommand{\dd}{\mathrm d}

\newcommand{\ie}{\emph{i.e.\ }}

\title{Randomized polynuclear growth with a columnar defect}

\author{Vincent Beffara \and Vladas Sidoravicius \and Maria Eulalia Vares}

\institute{V. Beffara\at
UMPA -- ENS Lyon,
              46 All\'ee d'Italie,
              F-69364 Lyon Cedex 07,
              France \\
              \email{vbeffara@ens-lyon.fr}
            \and
          V. Sidoravicius \at
              CWI, Kruislaan 413, 1098SJ Amsterdam, The Netherlands and
          \\
              IMPA,
              Estrada Dona Castorina 110,
              Rio de Janeiro RJ 22460-320,
              Brasil \\
              \email{v.sidoravicius@cwi.nl,  vladas@impa.br}
            \and
          M. E. Vares \at
              CBPF,
              Rua Dr. Xavier Sigaud 150,
              Rio de Janeiro RJ 22290-180,
              Brasil \\
              \email{eulalia@cbpf.br}
}

\date{Last modified: 03.02.09}

\begin{document}

\maketitle

\begin{abstract}
  We study a  variant of poly-nuclear growth where  the level boundaries
  perform  continuous-time, discrete-space random  walks, and  study how
  its  asymptotic behavior  is affected  by the  presence of  a columnar
  defect  on  the line.  We  prove that  there  is  a non-trivial  phase
  transition in the strength of the perturbation, above which the law of
  large numbers for the height function is modified.
\end{abstract}

\section{Preliminary considerations and statement of the results }


\medskip



Rigorous study of growth processes, and constant attempts to give
mathematical rigor to Kardar-Parisi-Zhang  theory and its
predictions, led to a very rich flow of results.  Remarkable, but
yet very limited progress was achieved in the last few years for
$1+1$ dimensional models by using a broad spectrum  of techniques
and arguments from  the theory of random matrices, first passage
percolation and interacting particle systems. An important role in
the successful application and interpretation of obtained results
was played by the fact that some properties of a particular 1+1
dimensional growth could be mapped and studied in terms of the
length of a maximal increasing sub-sequence in permutation of $N$
elements, with $ N \to \infty$, or as properties of the
one-dimensional totally asymmetric simple exclusion process, or
finally as a directed polymer in two dimensions subject to a random
potential.

\smallskip

\noindent The starting point and motivation of the present work was
the following question: how can a localized defect, especially if it
is small, affect the macroscopic behavior of growth system? This is
one of the fundamental questions in non-equilibrium  growth: is the
asymptotic  shape changed (\emph{faceted})  in the macroscopic
neighborhood of such a defect at any value of its strength, or, when
the defect is too weak, then  the fluctuations of the bulk evolution
become predominant and destroy the effects of the obstruction in such a
way that its presence becomes macroscopically undetectable?

\begin{figure}[h]
  \begin{center}
  \includegraphics[scale=.43]{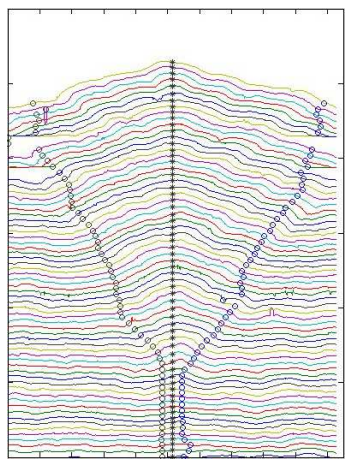}
  \end{center}
  \caption{Example: effect of the local defect on propagation of slow-combustion fronts: Higher concentration of potassium nitrate along
  vertical line changes shapes of otherwise nearly horizontal front
  lines. Higher velocity of propagation of combustion front along this
  line affects macroscopic picture producing a faceted shape around the
  defected region. (Experiment produced by M. Myllys et al.\ \cite{M-N}. Reproduced by permission.)}
\end{figure}

\noindent Such a vanishing presence of the macroscopic effect as a
function of the strength of obstruction represents what is called
\emph{dynamic phase transition}. The existence of such a transition,
its scaling properties, the shape of the density profile near the
obstruction, and also whether information percolates through the
obstruction, are among the most important issues. As we already
mentioned, in many  cases these questions can be translated into the
language  of pinning  of a directed  polymer in presence  of  bulk
disorder, or  to the question of whether the presence of a static
obstruction in a  driven flow, such as a  slow bond in one
dimensional totally asymmetric simple exclusion process (TASEP),
always results in change of the current, as predicted by mean-field
theory (see \cite{lebowitz:slowbond}).

\smallskip

\noindent  Similar questions are  often raised in the setting  of
first-passage percolation (FPP) (see \cite{kesten:fpp} for a
review). More specifically, consider an FPP  process in which the
passage times are shortened (in some random or deterministic way)
along one of the axes. Such type of perturbation usually is referred in physics
literature as
\emph{columnar defect}. Is it then true that the asymptotic shape of the
cluster is always  modified, or  is  there a phase transition for
the amount  of perturbation?  Fig. 1 shows experimentally obtained
successive slow-combustion fronts on the piece of paper, in the presence of the columnar defect (see \cite{M-N} for details).

\smallskip

\noindent All these questions remain rather poorly understood on the
experimental and the theoretical physics level, and essentially no
rigorous mathematical results are available.\footnote{In this context it also interesting to mention a
recent work (see \cite{Zeitouni}), where alternative methods to detect localized inhomogeneity were used. In \cite{Zeitouni} the following question was studied: Consider a graph with a set of vertices and oriented edges
connecting pairs of vertices. Each vertex is associated with a
random variable and these are assumed to be independent. In this
setting, the authors consider the following hypothesis  testing
problem: under the null, the random variables have common
distribution $N(0,1)$, while under the alternative, there is an
unknown path along which random variables have distribution
$N(\mu,1), \; \mu
>0$, and distribution $N(0,1)$ away from it. The authors study for
which values of the mean shift $\mu$ one can reliably detect the
presence of such defected path and for which values this is
impossible.} One particular model, however, known as
Poly-Nuclear Growth (or PNG for short) makes the exception. In this
process, on a flat substrate, layers appear on top of each other
according to a Poisson process, and then grow laterally with
constant speed $-1$ and $+1$, respectively, and merging whenever
they collide. This model is interesting in the sense that it lends
itself to the explicit computation of various quantities, such as
growth speed and interface fluctuations (see
\cite{johansson:png,spohn:png}).  The proofs depend on exact
algebraic computations, and it is not clear how to apply them as
soon as the dynamics are modified in a non-translation-invariant way
(with an exception of certain symmetrized cases (see
\cite{BR:symmetrized})). To treat the PNG case, new arguments have been
developed in \cite{beffara:BSSV} and \cite{beffara:PNG}.

\medskip

\noindent The  process which we introduce and study in this  paper,
the \emph{Randomized Poly-Nuclear Growth} (or RPNG for short), has
several apparent similarities with PNG, but the substantial difference is
that after the nucleation, boundaries of layers (plateaux) move as
independent random walks instead of moving deterministically. This
dynamics corresponds to the modified Glauber dynamics for the
two-dimensional Ising model  on a solid surface at temperature $0$,
as considered in \cite{chayes:css},  with additional nucleations.
Informally it can be described  as follows:  On each  edge  of the
discrete  line $\mathbb  Z$, with  some positive rate $\lambda$,
pedestals of width and height $1$ appear. Every pedestal  will
evolve  into   a  plateau,  its   endpoints  performing independent
continuous-time, discrete-space symmetric random walks with jump rate $1$,
until one of the following happens:
\begin{itemize}
\item If two  plateaux \emph{collide}, \emph{i.e.} if the  right edge of
  one occupies the same site as the left edge of the next one, then they
  merge and only  the two extremal endpoints continue  to perform random
  walks;
\item If one of the two endpoints  of a plateau jumps to the location of
  the  other one,  so that  the plateau  shrinks to  width $0$,  then it
  vanishes.
\end{itemize}

If a pedestal appears on an edge  which is already spanned by one of the
plateaux, it  appears on  an additional  layer on top  of it.   Only two
plateaux  evolving within  the same  layer  can merge,  and besides  the
growth of  the plateaux  at different layers  has to obey  the following
constraint: Every plateau  which is not on the  ground level must always
lie entirely on top of a plateau at the layer below it.  Each jump which
would break  this constraint is simply suppressed. For  any fixed value
of  $\lambda$  (see Proposition~\ref{prp:lln}  below) the  process has linear
growth, \emph{i.e.}   if by $h_{t}(e)$ we denote  the height of the
aggregate at  time $t$ above the edge $e$, then the following limit
exists:
\begin{equation}
  v(\lambda) := \lim_{t \to +\infty} h_{t}(e)/t >0.
\end{equation}
Notice that in the case  of a strictly growing  process (such as first-
or last-passage percolation), the law of large numbers is usually a
direct consequence  of a  statement of  sub-additivity, while  if the
dynamics is allowed to both grow  and shrink, like in the present case,
even the existence of the asymptotic speed is non-trivial.


Next we  will introduce  a \emph{``columnar defect''}  to the system
by changing the nucleation  rate on the edge  $e_0 := <0,1>$  from
$\lambda$ to $\lambda+\lambda_0$, with $\lambda_0>0$, and  the main
question which we address is how does this affect the asymptotic
growth rate  of the aggregate  above the origin. More precisely, let
$h_{t}^{\lambda_0}(e_0)$ be the height  of the aggregate above $e_0$
at  time $t$  under  the modified dynamics (we make the dependency
on $\lambda$ implicit in this notation, since the parameter of
interest  is $\lambda_0$). Assuming that the following limit exists
\begin{equation}
v(\lambda,\lambda_0) := \lim_{t  \to +\infty}
h_{t}^{\lambda_0}(e_0)/t,
\end{equation}
we are interested in whether the inequality
\begin{equation}
v(\lambda,\lambda_0)  > v(\lambda) \equiv v(\lambda,0)
\end{equation}
holds for all $\lambda_0 > 0$ and set:
$\lambda_c (\lambda) :=  \inf \{  \lambda_0 :  v(\lambda,\lambda_0) >
    v(\lambda) \}$ be the critical value of $\lambda_0$.

We  are  now ready  to state  our main result.

\medskip

\begin{thm}\label{thm:main}

$$
    v(\lambda,\lambda_0)    :=     \lim_{t    \to    +\infty}    \frac
    {h_{t}^{\lambda_0}(e_0)}  t=\lambda + \max\{\lambda_0 -1,0\}$$
    almost  surely. In particular,
    $\lambda_c (\lambda) = 1$, for any value of $\lambda$.

\end{thm}

In proving this  theorem, we actually obtain a  more precise
description of the  behavior of the system  in the supercritical
phase $\lambda_0 > \lambda_c (\lambda)$: the behavior  of
$h_{t}^{\lambda_0}(e)$  for $e \neq e_0$ is perturbed only
sub-linearly.  In other words, the defect creates a kind of
``antenna''  at the origin. This contrasts with the expected
behavior  of the usual PNG with  a supercritical localized defect,
for which the region of influence of the additional nucleation
process is believed to extend linearly in space. However if in the RPNG model
we would remove assumption of the symmetry of walks and would take the
boundaries of plateaux moving as random walks with a drift (left boundary
with a negative drift and right boundary with positive one), we expect
the behavior of the system to resemble the behavior of PNG.

Moreover, the phase transition is non-trivial even without
nucleation in the bulk (\emph{i.e.} $\lambda = 0$).  Even in this
simplified situation, for any value $\lambda_0 > 0$ the origin is
eventually covered  by plateaux for  all times, however $0 <
\lambda'_c := \lambda_c (0) < + \infty $  so that if $\lambda_0  <
\lambda'_c $ then $v(0,\lambda_0)=0$, \emph{i.e.} there is
sub-linear growth,  while for $\lambda_0 > \lambda'_c $  we have
$v(0,\lambda_0)>0$. As we shall see, $\lambda'_c=1$, in agreement
with Theorem 1.

\medskip

\noindent {\bf Comment} It might be instructive to compare with the
following simplified model in the half-line $\mathbb R_+$. Consider
a situation corresponding to $\lambda=0$ and an additional vertical
wall at vertex $0$ to which the left endpoint of each pedestal is
fixed: the initial configuration is empty; at each occurrence of a
Poisson process  in $\mathbb R_+$ of intensity $\lambda_0$, a
new pedestal is created at the edge $e_0$ --- and this is the only
source in the system. Right boundaries of pedestals move according
to the same mechanism as described in the basic homogeneous  model,
except that now each pedestal is annihilated upon its right boundary
hitting the site $0$.

Of particular interest is the particular case $\lambda_0=1/2$, when
the system is in exact bijection with $0$-temperature Glauber
dynamics for the two-dimensional Ising model, with an initial
configuration consisting in the first quadrant uniformly at spin
$-1$ and the rest of the lattice uniformly at spin $+1$. The height
function at $e_0$, here represented by $N_t$, the number of
pedestals in the system at time $t$,  describes the shape of the
$+/-$ interface, and in the large time limit it is known to behave
according to Lifshits' law of motion by mean-curvature (see
\cite{chayes:css}, \cite{spohn:motion}).

In this particular case, there is also a one-to-one correspondence
with the simple symmetric exclusion process:
$N_t$ corresponds to the position of the rightmost particle in a
symmetric simple exclusion process in $\mathbb Z$ starting with the
configuration where the occupied sites are exactly those $x \le 0$.
A classical result of Arratia (see Theorem 2 in \cite{arratia:arr1},
cf.\ also \cite{liggett:book1}, comments at p.\ 414) tells in particular that
$N_t/t \to 0$ almost surely. A simple comparison with this case tells us
that $1/2$ is the critical value of $\lambda_0$ in this simplified
model.

\noindent{\bf Structure of the paper.}  In
Section~\ref{sec:construction} we give  a formal construction  of  the
process,  in terms of a Harris-type representation, and discuss two
other equivalent representations. In Section~\ref{sec:hom} we derive
basic properties of the homogeneous process without the defect at
$e_0$. More specifically,  we prove a law of large numbers for the
height of the  interface above edge  $e_0$; a key  step towards this
result is the  proof  of  eventual  covering of  any finite-size
interval of the line. Finally, in Section~\ref{sec:full} we present
the proof of Theorem \ref{thm:main}.

\section{Construction of the process: three ways to look at the system}
\label{sec:construction}
\subsection{Description as an interacting particle system}

One  can describe this  process  in terms  of  an interacting
particle system.  Each left endpoint of  a plateau is represented by
a \emph{left-particle} or $\ell$--particle for short, at the
corresponding site of $\mathbb Z$, and each right endpoint of a
plateau as a \emph{right-particle}, or $r$--particle (which can be
interpreted as an anti-particle). If several plateaux have their
endpoints at the same site,  this site will contain several
particles, but notice that a site will never contain particles of
both types.

The dynamics then looks like  a zero-range process with annihilation
and nucleation, \emph{i.e.}:
\begin{itemize}
\item For each non-empty site, with rate $1$, one of the particles jumps
  with equal probability  to a neighboring site; and  if the target site
  contains a particle  of the other type, then  the jumping particle and
  one of the target particles annihilate instantaneously;
\item On every edge $(x,x+1)$ of $\mathbb Z$, with rate $\lambda$, there
  is nucleation, meaning that an $\ell$--particle appears at $x$ and
  an $r$--particle appears  at $x+1$; in each case, the particle is
  instantaneously  annihilated (together with a particle of the opposite type),
  if a particle of the opposite type is already present at the corresponding site.
\end{itemize}


\noindent\textbf{Remark.} Notice that the particle
configuration at time $t$ gives only partial information about the
initial system; in particular, it gives no way of following the
evolution of a given plateau boundary as a function of time --- indeed,
with our description of the particle dynamics, the notion of a tagged
particle is ill-defined. Nevertheless, knowing the particle
configuration at time $t$ allows for the reconstruction of the plateaux
configuration at the same instant up to a global height change; and of
course, the full space-time history contains exactly the same
information for both descriptions of the system.

\medskip

\noindent We first construct the above described  interacting
particle system  using  \emph{Harris graphical representation}. Let
$\EE:=\{<x,x+1>, x\in\ZZ\}$ be the  set of edges of the integer line
$\ZZ$.  Define the following three independent Poisson systems:
\begin{itemize}
\item $\LL_r=\{\LL_r (e),  e\in\EE\}$ is  a collection of  independent Poisson
  processes  of   left-arrows,  each  with  intensity   $\dd  t/2$;  the
  individual arrows (or \emph{marks}) of the process at edge $e$ will be
  denoted by $\LL_{r}^n (e)$, $n\geqslant 1$, with similar notations for the
  other two systems;
\item $\LL_{\ell} =\{\LL_{\ell}{(e)}, e\in\EE\}$ is  a collection of independent Poisson
  processes of left-arrows, each with intensity $\dd t/2$;
\item $\RR_r=\{\RR_r{(e)}, e\in\EE\}$ is  a collection of independent Poisson
  processes of right-arrows, each with intensity $\dd t/2$;
\item $\RR_{\ell} =\{\RR_{\ell}{(e)}, e\in\EE\}$ is  a collection of independent Poisson
  processes of right-arrows, each with intensity $\dd t/2$;
\item $\NN=\{\NN_{e}, e\in\EE\}$ is  a collection of independent Poisson
  processes of undirected edges, each with intensity $\lambda \dd t$.
\end{itemize}
$\LL_r$ (respectively $\LL_{\ell} $) will correspond to jumps of
right particles (respectively left particles) to the left, $\RR_r$
and $\RR_{\ell}$ to jumps of right and, respectively, left particles
to the right, and $\NN$ to nucleation times.

Given $T>0$, we first construct the  process up to
time $T$,  as follows.  There is positive  density of edges  carrying no
mark of any of the Poisson  before time $T$; between any two such edges,
there are almost surely only  finitely many marks appearing by time $T$,
so  that they  can be  ordered  chronologically, and  we  only need to
describe the action of each of the marks.

\begin{itemize}
\item If there is  a directed arrow of type $r$ or, respectively, type $\ell$
at edge $e$ and  time $t$, and if at this time the site  at the
source of the arrow is not empty and contains at least one right or,
respectively, at least one left particle, then one particle at this
site jumps in the direction  of the arrow; moreover,
  if after  this jump  it meets  a particle of  the opposite  type, they
  annihilate instantaneously.
\item If there is a mark of $\NN$ at edge $e=<x,x+1>$ and time $t$, then
  at this time an $\ell$--particle appears at site $x$ and an $r$--particle
  appears  at site  $x+1$; in each case if the  corresponding site contain
  particles of  the  opposite  type,  a  pair of opposite type particles
  is  annihilated instantaneously.
\end{itemize}

It is then easy to see  that the obtained processes for different
values of $T$ are consistent, and letting  $T$ go to infinity gives
the process for all time.  We will refer to this version of  the
construction as the process in  \emph{quenched environment},
\emph{i.e.}  seeing the Poisson processes as a random scenery in
which a deterministic process evolves.

\noindent {\bf Remark.} Just for the sake of constructing of the
process, there is no need to distinguish between $\LL_r$ and
$\LL_\ell$, or $\RR_r$  and $\RR_\ell$. The distinction plays a role
only when describing the evolution of the symmetric difference between two
configurations, which we shall need later.

\medskip

\noindent\textbf{Particle system corresponding to the process with a
columnar defect:}
A columnar defect in the plateaux dynamics, discussed in the Introduction, is described
in the interacting particle system formalism  by the addition  an auxiliary nucleation
process $\NN_0$  with rate $\lambda_0>0$ localized on the edge
$e_0$. In other words, the pedestals appear with rate $\lambda$
everywhere except at $e_0$ where the nucleation rate is $\lambda +
\lambda_0$.

\noindent The construction of the new process  can be done in exactly the same
way as for the  unperturbed system, using a graphical construction
and letting the marks of $\NN_0$ act the same way as those of $\NN$.


\noindent \textbf{Notation:} $n_x(t)$ will denote the number of
particles at site $x$ and time $t$ in the above construction;
$\NN_{e}(t)$ denotes the number of nucleation marks occurring at the
edge $e$ up to time $t$.

\subsection{Construction in terms of interacting random walks}
\label{sub:1.2}

In this subsection we will give an alternative
interacting-random-walks-type description of our system. To do so,
one needs a rule to  determine how do these walks interact, and
which of the particles at a given  site  will  jump  when  an  arrow
(in the previous construction) appears,  and  which  will  be
annihilated when a  particle of the opposite type  arises. Obviously
the process  itself will  not be  affected by  the choice,  so we
have some freedom here. Some  of the more natural choices are  to
declare that the oldest  particle (\emph{i.e.}  the  one that
appeared  earliest in  the system), or the youngest particle, or the
one that has been at the given site for the  longest time, or the
shortest time,  is ``active". The rule can even be different for
jumps and annihilations. We will refer to such a rule as a
\emph{priority rule}.


\noindent In particular, the previous graphical  construction can be
replaced  with the following interpretation  of the dynamics:

\smallskip

\noindent \emph{1) Evolution of a single particle.}  Whenever a
particle  $i$ is born into the system, it  comes together with a
``project", \emph{i.e.} a family  $(T_i^n)$ of  independent
exponential  times of  rate $1$  and a family $(\varepsilon_i^n)$ of
independent signs, and a clock. Those will represent waiting  times
between consecutive jumps, and, respectively,  directions of jumps
for the path  of this particle, if it were moving alone in an
otherwise empty space.
\smallskip

\noindent \emph{2) Interaction.}  Choose any priority rule. When a
particle is alone, it is always ``active'', \emph{i.e.} its clock
runs deterministically with speed $1$, and it executes its
project, i.e.\ jumps according to times and directions prescribed in
its project. When there are several particles at a given site, only
one of them is active, according to the chosen rule, and the others
are ``sleeping''. The clock of the active particle runs
deterministically with speed $1$, the  clocks of  the other
particles do not run at all.  Whenever the clock of any particle
reaches the next  waiting time for the project  of this particle,
the particle jumps and the clock is reset.

\smallskip

\noindent {\bf Remark.} Another natural choice  for a jumping rule
could be that  at a site with $n$ particles, the clock of each of
them runs with speed $1/n$; again the distribution of the numbers of
particles at different  sites would be the same as in the graphical
construction.
\smallskip

\noindent  The formal construction of this dynamics could be
done similarly as in the Harris construction: Given a space-time
configuration of particles (space time Poisson process,
corresponding to the $\NN$ system in the graphical construction) and
their families of ``projects", for any fixed $T>0$, there are
infinitely many edges to the left and to the right of the origin
such that by time $T$ they will not be crossed by any type of
particle, neither a nucleation will occur at them by time $T$. It
follows quite simply that the evolution of the system within the
space time block in between two such edges and $[0,T]$ depends only
on finitely many particles, being independent from the evolution
outside.




\subsection{Back to the PNG-like description}

We saw how to go from  a PNG-like description to a family of
interacting particles, but our  final goal is to describe the long
term behavior of the height  (\emph{i.e.} the number  of layers) at
a given edge,  so we need a  way to go  back to  the layer
structure.  As we saw  already the $\ell$--particles correspond to
left  boundaries  of  plateaux, and the $r$--particles correspond to
right boundaries. It only remains to compute the height at  edge
$e_0=<0,1>$.

Let $(\phi_t^\ell)$ be  the flow of $\ell$--particles through the
edge $e_0$:  $\phi_0^\ell=0$ and  $\phi^\ell$ increases (resp.\
decreases)  by $1$ whenever an $\ell$--particle crosses $e_0$ to the
right (resp.\ to the left).  Define $(\phi_t^r)$ similarly as the
flow of $r$--particles. Let $h_t(e_0)$ be the height of the  profile
at edge $e_0$ and time $t+$ (\emph{i.e.}, the càdlàg version in time
of the height process). Then, recalling previous notation,
\begin{equation}
  h_t(e_0) = \NN_{e_0}(t) - \phi_t^\ell + \phi_t^r,
  \label{htatzero}
\end{equation}
From $h_t(e_0)$ and  the previous discussion, it is  then easy to
obtain the whole profile of plateaux.

\noindent {\bf Remark.} The graphical method gives us a construction
of $h_t := \{ h_t (e)\}_{e \in \EE}$ as a Markov process. The
$\{h_t\}$--evolution is attractive (or monotone in the sense of
definition 2.3 of Chap. II in \cite{liggett:book1}), though the
particle system is not. The height process is also monotone in the
nucleation rates.

\medskip


\section[The homogeneous system]{Behavior of the homogeneous process ($\lambda > 0, \; \lambda_0 = 0$)}
\label{sec:hom}

\subsection{Level lines, Hammersley--process--like description}

Let $L_t$ (resp.\ $R_t$) be the  position of the left (resp.\
right) end of the first-layer plateau containing $e_0$ at time
$t$, if there is one \emph{i.e.} if  $h_t(e_0)\geqslant 1$ --- for
completeness,  one can let $L_t=+\infty$ and $R_t=-\infty$ if
$h_t(e_0)=0$.

Let $\tau$  be the first  time at which  $R_t=L_t+2$, if any
(otherwise $\tau=+\infty$).  In the absence of other layers
and of further nucleation, the behavior after $\tau$ of $L_t$ and
$R_t$ until the first time when $R_t=L_t+1$ would exactly be that of
two independent, centered continuous-time random walks. The presence
of the second and above layers, by the previous description, can
only affect the pair by the suppression of left jumps of $R_t$ and
right jumps of $L_t$
--- notice that if, in the description at Subsection 2.1, there are several
particles at  $R_t$ ($L_t$)  at time $t$ and one has to jump to
the right (left, resp.), then so does $R_t$ ($L_t$ resp.).

Moreover,  with   rate  $\lambda$  there  is  nucleation   on  the
edge $<R_t,R_t+1>$,  and when this  occurs, $R_t$  jumps to  the
right  by at least one step.  The last thing that can affect $R_t$
is if the plateau $[L_t,R_t]$ merges  to the  next plateau to  its
right, which  again can only increase the  value of $R_t$. The
same remarks apply \emph{mutatis mutandis} to $L_t$.

Taking  all those  effects into  account, one  obtains that  the
process $(R_t)$ stochastically  dominates the biased random  walk
$(\tilde R_t)$ that has  rates $1/2$ to jump to  the left and
$1/2+\lambda$  to jump to the right.  Similarly, $(L_t)$  is
stochastically dominated  by $(\tilde L_t)$ which is biased to the
left. Since the probability
\begin{equation}
 p=P(\forall t>0,\; \tilde R_t - \tilde L_t > 1)
  \label{eq:eventual}
\end{equation}
is positive, using standard arguments we see that for $t_0$ large enough we have
\begin{equation}
  P \left( \forall t \geqslant t_0, \; R_t > \frac {\lambda t}{2} \;
  \mathrm{and} \; L_t < - \frac {\lambda t}{2} \right) > 1- \const e^{-\const t_0},
\end{equation}
for suitable positive constants $c_1, c_2$. Applying the
Borel-Cantelli Lemma, one in particular sees that a.s.\ any compact interval of
the real line will eventually remain covered  by the plateau $[L_t,R_t]$.
Moreover, extending the same line of reasoning one sees that
$$
 \liminf_{t \to \infty} \frac 1 t\log R_t \ge C \quad a.s.
$$
for suitable constant $C>0$.



\bigskip

\noindent We  shall give  an additional  description  of the  time
evolution  of the  system, in  terms of  a family  of  space-time curves
similar to a Hammersley process (cf.\ \cite{AD}). The main difference is
that the paths of the  particles between meeting points are random walks
instead of straight lines.

More precisely,  one can see  the path of  any given particle,  from
its birth  to its  eventual  annihilation,  as a  sequence  of
vertical  and horizontal line  segments. Let  $\CC_1$ be the  union
of  the space-time paths of all the particles corresponding to
endpoints of plateaux of the first layer, together with the segments
joining the birth places in space-time of the corresponding pairs.

\begin{figure}[h]
  \begin{center}
  \includegraphics[scale=.8]{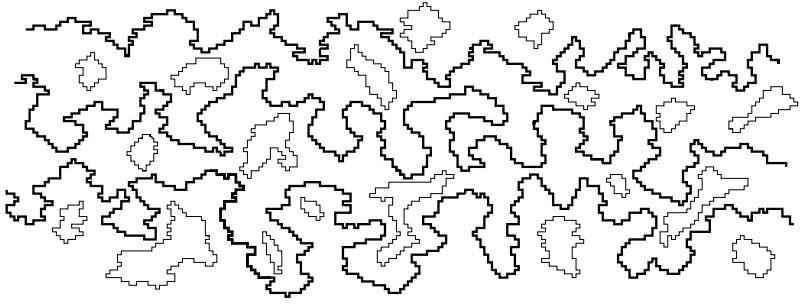}
  \end{center}
  \caption{Illustration of an RPNG process. (Level lines are bold.)}
  \label{fig:rpng}
\end{figure}


As we have seen above, for each bounded segment $[-K,K]$ there is an a.s.\ finite random time
$t_K$, such that, the segment $[-K,K]$ will be covered by a plateau of the first level
for all $t\ge t_K$. This means that $\CC_1$
is  the   disjoint  union  of   a  countable  family  of   finite  loops
(corresponding to  plateaux with a  finite life span) and  one unbounded
curve,  which we will  denote by  $C_1$ and  call the  \emph{first level
  line} of the process.

Similarly, let $\CC_k$  be the union of all the  space-time paths of
the endpoints of plateaux  of the $k$th layer together with the
segments joining the birth places in space-time of the corresponding
pairs.  $\CC_k$ is then the disjoint union of a  countable family of
finite loops, and one unbounded curve $C_k$. Notice that in general
a component of $\CC_k$ and a component of $\CC_l$ for $l\neq k$ need
not be disjoint.

\bigskip

Using this description, it is possible to compute the height of the
RPNG process at any  edge $e$ and (almost) any time $t$:  Indeed, it
is equal to the number of level lines that  separate $(e,t)$ from
the line $t=0$ (\emph{i.e.} that pass below it an odd number of
times), plus the number  of loops containing  $(e,t)$. This means
that the law  a of large numbers  for $h_t$  can be translated  into
a statement  about the asymptotic density of the level lines.

\subsection{The asymptotic speed}

The first step in proving a strong law of large numbers for the
homogeneous process is to obtain an \emph{a priori} bound for the
number of particles present at a given site, at a given time. We
obtain such a bound in the form of a stochastic domination result.


\begin{lem}
  \label{lem:comprw}
  There exist three countable families $(X^x_{\ell,\cdot})_{x\in\ZZ}$, $(X^x_{r,\cdot})_{x\in\ZZ}$ and
  $(Y^x_\cdot)_{x\in\ZZ}$ of integer-valued stochastic processes, and a
  coupling of them with the homogeneous particle system introduced above
  satisfying the following conditions:
  \begin{enumerate}
  \item For each $x\in\ZZ$, $(X^x_{\ell, t})_{t\geqslant0}$ and $(X^x_{r, t})_{t\geqslant0}$ are simple symmetric random
    walks on $\ZZ_+$, with reflection at the origin, jump rate $1$ at the origin and $2$ elswhere,
    and $(Y^x_t)_{t\geqslant0}$ is a simple symmetric random walk on $\ZZ$ with
    jump rate $\lambda$;
  \item Almost surely, for all $x\in\ZZ$ and $t\geqslant0$,
    \begin{equation}
      n_x(t) \leqslant X^x_{\ell, t} + X^x_{r, t} + \sup_{s,s'\in[0,t]} (Y^x_s - Y^x_{s'}).
      \label{27July}
    \end{equation}
  \end{enumerate}
\end{lem}

\begin{proof}
  Let $x\in\ZZ$, and in the  graphical construction of the process (cf.\
  Section \ref{sec:construction}),    let    $(T_{\ell, x,k})_{k>0}$   (resp.\
  $(T'_{\ell, x,k})$) be the  $k$th time at which there  is an $\ell$-arrow pointing
  toward  (resp.\  away  from)  site  $x$.   Define  $X^x_\ell$  as  follows:
  $X^x_{\ell, 0}=0$, $X^x_\ell$  jumps up  by $1$ at  each  $T_{\ell, x,k}$, and
  jumps down  by $1$ at each $t\in\{T'_{\ell, x,k}\}$  such that $X^x_{\ell, t-}>0$.
  Clearly, $X^x_\ell$ defines a continuous-time random walk as in the statement. Similarly
  we define  $X^x_r$ using $r$-arrows. Let $Y^x_t = \NN_{\langle x,x+1 \rangle} (t)- \NN_{\langle x-1,x
   \rangle} (t)$.

   $n_x$ can  change for two reasons (nucleation or jumps): when a
   particle at $x$ jumps away, $n_x$ decreases by $1$, and  when a
   particle  at $x\pm1$ jumps onto $x$, it either increases by 1, or (if
   $n_x>0$) it decreases by $1$ depending on the types of particles
   involved;  if there is nucleation at the edge $\langle x-1,x\rangle$
   or $\langle x,x+1\rangle$, then $n_x$ again increases or decreases by
   one according to the types of particles involved, as before.

    Fix $t>0$, and let $$\tau := \sup \{ s<t : n_x(s)=0
   \}.$$  If $\tau=t$ and $n_x(t)=0$ we have nothing to prove. From
   now on we are assuming that this is not the case. In particular, during
   the time interval $[\tau, t]$, the site $x$ is never empty and the type of particles
   present in it does not change; besides, $n_x(\tau-)=0$.

   Apply the particle-counting of the previous paragraph between times
   $\tau$ and $t$. The contribution from nucleations is exactly equal to
   $\pm(Y^x_{t} - Y^x_{\tau-})$ (where the sign depends on the type of
   particles present at site $x$), and the contribution from jumps is
   controlled by the jumps of $X^x_{\ell, t}$ or $X^x_{r, t}$. This leads to the
   following estimate: $$n_x(t) \leqslant |Y^x_t - Y^x_{\tau-}| + X^x_{\ell, t} + X^x_{r, t}.$$
   The  bound in the statement of the  Lemma directly
   follows.
\end{proof}

%


\begin{prp}
  \label{prp:lln}
  Recall that $e_0=<0,1>$. For every  $\lambda>0$,
  almost surely,
  $$\lim_{t\to \infty}\frac {h_t(e_0)} {t} = \lambda.$$
\end{prp}

\begin{proof}
  Let $Q>0$ be a positive integer, and let
  $$H_t(Q) := \sum_{x=-Q}^{Q-1} h_t(<x,x+1>)$$
  be the \emph{integrated height} of the process over the interval
  $[-Q,Q]$. The variations of $H_t(Q)$ in time come from several sources
  (according to the definition of the height process):
  \begin{itemize}
    \item Whenever there is a nucleation on one of the edges in
      $[-Q,Q]$, $H_t(Q)$ increases by $1$;
    \item Whenever an $\ell$--particle  in $[-Q+1,Q-1]$
      jumps to the right (resp.\ to the left), $H_t(Q)$ decreases (resp.\ increases)
      by $1$
      \item Whenever an $r$--particle in $[-Q+1,Q-1]$
      jumps to the left (resp.\ to the right), $H_t(Q)$ decreases (resp.\ increases)
       by $1$;
     \item Whenever a particle sitting at $Q$ (resp.\ $-Q$) jumps to site
   $Q-1$ (resp. $-Q+1$), $H_t(Q)$ changes by $+1$ or $-1$ according to the
type of particle ($\ell$ or $r$) involved.
\end{itemize}
  Any other jump or nucleation might change the value of $n_{\pm Q}(t)$,
  but it will never affect the value of $H_t$. Notice in the second and third
  cases, the decrease of $H_t(Q)$ occurs independently of the jump leading
  to an annihilation or not.

  Notice that a particle sitting in $[-Q+1,Q-1]$ jumps to the left and
  to the right with the same rate (be it an $\ell$- or an $r$-particle),
  and that the total jump rate of
  particles in this interval is, by construction, bounded by $2Q-1$ (the
  total jump rate of particles at any given site being $1$ if the site
  is non-empty, and of course $0$ if the site is empty). Hence, the
  contribution of such jumps to $H_t(Q)$ is at most of the order of the
  position at time $t$ of a symmetric, continuous-time random walk
  jumping with rate $2Q-1$ --- \ie, they are of order $\sqrt t$.

  What remains is the contribution of nucleations inside $[-Q,Q]$,
  which have total rate $2Q\lambda$ and satisfy a strong law of large
  numbers, and the fourth term which, although difficult to estimate
  precisely due to its dependency on the exact particle distribution,
  only contributes to $H_t(Q)$ with rate at most $1$ (\ie, the sum of
  the jump rates from $-Q$ to $-Q+1$ and from $Q$ to $Q-1$).

  Summing all the above contributions to $H_t(Q)$, and then letting $t$
  go to $+\infty$, we obtain with probability $1$:
  \begin{equation}
    2Q\lambda-1 \leqslant \liminf_{t\to\infty}\frac{H_t(Q)}{t} \leqslant \limsup_{t\to\infty}
    \frac{H_t(Q)}{t} \leqslant 2Q\lambda+1.
    \label{eq:liminfsup}
  \end{equation}

  \medskip

  Now, we know that at time $t$, the difference in height between two
  neighboring edges, say $<x-1,x>$ and $<x,x+1>$, only depends on the
  number (and type) of particles present at site $x$ at time $t$. By (\ref{27July})
  we see that almost surely, for every $x\in\ZZ$,
  $$\lim_{t\to\infty} \frac {n_x(t)} {t^{}} = 0.$$

  The height $h_t(e)$ at any edge $e\in[-Q,Q]$ can be computed from the height
  at edge $e_0$ by only looking at the values of the $n_x(t)$ for
  $x\in[-Q+1,Q-1]$. In particular, for such an $e$,
  $$\left| h_t(e) - h_t(e_0) \right| \leqslant \sum_{x=-Q+1}^{Q-1}
  |n_x(t)|.$$
  Summing over $e$ to obtain $H_t(Q)$, and taking a crude upper bound,
  leads to:
  $$\left| H_t(Q) - 2Q h_t(e_0) \right| \leqslant 2Q \sum_{x=-Q+1}^{Q-1}
  |n_x(t)|.$$
  Combining this with equation (\ref{eq:liminfsup}) and the previous
  estimate on the $n_x(t)$, we obtain, still with probability $1$:
  $$2Q\lambda-1 \leqslant \liminf_{t\to\infty}\frac{2Q h_t(e_0)}{t} \leqslant \limsup_{t\to\infty}
    \frac{2Q h_t(e_0)}{t} \leqslant 2Q\lambda+1$$
  (notice that we are still working with fixed $Q$). Dividing by $2Q$,
  we get
  $$\lambda-\frac 1 {2Q} \leqslant \liminf_{t\to\infty}\frac{h_t(e_0)}{t} \leqslant
  \limsup_{t\to\infty}\frac{h_t(e_0)}{t} \leqslant \lambda+\frac 1 {2Q}.$$
  Since this holds for every positive $Q$, the proposition is
  proven.
\end{proof}

\section{Proof of Theorem 1}
\label{sec:full}

We
start with a homogeneous RPNG process as described in Section~\ref{sec:construction}; let $h_t$ be its height function. Notice that
the whole dynamics can be described in terms of this function.
\noindent \textbf{Notation:} $e_x=\langle x,x+1\rangle$, for each integer $x$.

Now add an extra nucleation with rate $\lambda_0$ at the edge $e_0$  to the Harris graph
of Section 2, and  let $h'_t$  be the  height function  of the  modified  process, with
$n'_x$ being the number of particles at $x$ in the corresponding new particle system.  By
monotonicity of the plateaux dynamics (cf. Remark right after equation (\ref{htatzero})),
we have that  $h'_t \geqslant h_t$; and
let $\tilde h_t = h'_t - h_t \geqslant 0$ be \emph{the discrepancy process}. We are interested in
differences in the law of large  numbers between the two versions of the
process,  \emph{i.e.} in the asymptotic  properties  of
$\tilde h_t/t$.

The discrepancy process can be thought as a particle system, composed of what we call
\emph{virtual particles}: According to the difference
$$\tilde h_t(e_{x-1}) - \tilde h_t(e_{x})$$ being positive or negative
we say that at time $t$ there are  $$\tilde n_x(t) :=
|\tilde h_t(e_{x-1}) - \tilde h_t(e_{x})|$$ virtual right or, respectively, virtual left particles at
site $x$. In this interpretation the creation of virtual particles corresponds to a nucleation in ${\cal N}'_{e_0}\setminus {\cal N}_{e_0}$.

\bigskip

\begin{figure}[ht!]
  \begin{center}
    \includegraphics[scale=.3]{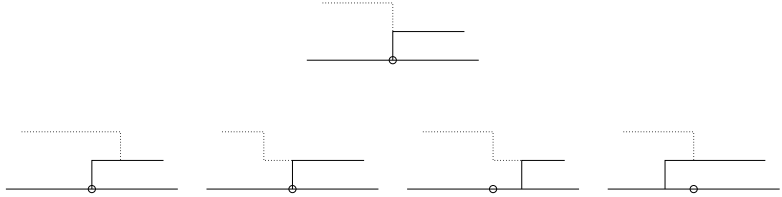}
  \end{center}
  \caption{Comparison of the perturbed and unperturbed dynamics in the
  case $b>a$. In plain lines are the height functions for the
  unperturbed system, and in dotted lines are the height functions for
  the perturbed one. On top is the initial configuration, and below are
  configurations after occurrence of an arrow in $\mathcal R_r$,
  $\mathcal L_r$, $\mathcal R_\ell$, and $\mathcal L_\ell$,
  respectively. The circle on the axis indicates the location of site
  $x$.}
  \label{fig:comparison}
\end{figure}

Using the Harris description of subsection 2.1 (or the random walk
representation of subsection 2.2) one can check explicitly that the
spatial displacements of virtual particles remain as of independent
random walks even in presence of the interaction with other virtual
particles while jump times can be possibly delayed --- or in other
words, that conditionally on the whole particle system, the current
rates of jump of a given virtual particle to the left or to the right
always remain equal.  Doing so in full detail involves looking at many
different cases, but the discussion can be shortened considerably by the
following considerations:
\begin{itemize}
  \item If a site contains only virtual particles, then it is obvious
    from the construction that the first jump of one of them away from
    the site occurs with equal rate to the left and to the right; so we
    only need to check what happens when normal and virtual particles
    share a site.
  \item If a site contains both normal and virtual $\ell$-particles,
    then the perturbed system also contains $\ell$-particles. In
    particular, the original and perturbed system will use the same
    arrows in the Harris construction, meaning that their difference
    will not be affected by the dynamics and the virtual particles will
    not move as long as there are normal particles present, at which
    point we are back to the previous case. Of course, this also applies
    to the case of $r$-particles.
  \item If a site $x$ contains $a>0$ normal $\ell$-particles and $b>0$
    virtual $r$-particles, then one needs to check $3$ sub-cases ($a<b$,
    $a=b$ and $a>b$), leading to different behavior in the perturbed
    system.  We treat one of them fully, the others being similar:
    Assume that $b>a$ (see Fig.~\ref{fig:comparison}), so that the
    perturbed system contains $(b-a)$ $r$-particles at site $x$. The
    occurrence of an arrow from $x$ affecting left-particles will affect
    the unperturbed system and not the perturbed one, meaning that one
    of the virtual particles will follow the jump as well; the
    occurrence of an arrow affecting right-particles will only modify
    the perturbed system, and again in this case a virtual particle will
    follow; in both cases, the symmetry of the arrow processes
    translates into symmetry for the jumps of virtual particles.
\end{itemize}

\bigskip

\noindent {\bf Proof of the Theorem.}
Due to the monotonicity, the height of
the process with additional nucleation at the origin stochastically
dominates that of the homogeneous process. In particular, we already know that
almost surely
$$\liminf \frac{h'_t(e_0)}t \geqslant \lambda.$$

\noindent {\bf Step 1.} Define
$$
\Delta_t =  h'_t(e_0) - \max\{ h'_t(e_{-1}), h'_t(e_1) \}.
$$
If $\lambda_0 > 3 + \lambda$, then
\begin{equation}
\label{liminf2}
\liminf_{t \to \infty} \frac{\Delta_t}{t} > 0.
\end{equation}
Indeed, if  $\Delta_t > 0$, the following can happen:
\begin{itemize}
\item if $h'_t(e_{-1}) = h'_t(e_1) $, then $\Delta_t$ may increase by one only because of a
nucleation at the edge $e_0$, which happens at rate $\lambda + \lambda_0$, and it decreases by one either when some particle jumps away from vertex $0$ or $1$, which happens
with rate $2$, or when nucleation occurs at one of edges $e_{-1}$ or $e_1$, which happens
with rate $2\lambda$, or if  a particle jumps from vertex $-1$ (resp.\
from vertex $2$) to $0$ (resp.\ to 1), which happens at rate 1 (if such a particle is present);
\item if $h'_t(e_{-1}) < h'_t(e_1) $, then $\Delta_t$ may increase by one either because of a
nucleation at the edge $e_0$, which happens at rate $\lambda + \lambda_0$, or because of a jump of a
particle from vertex $2$ to vertex $1$, provided a particle exists at vertex $2$, and which happens at rate $1/2$. On the other hand $\Delta_t$ decreases by one either if a particle jumps away from vertex $1$, which happens at rate $1$, or if a particle jumps from vertex $0$ to $1$, which happens at rate $1/2$, or, finally, if nucleation occurs at the edge $e_1$, which happens at rate $\lambda$.
\end{itemize}

\noindent Examining similarly the cases when $\Delta_t = 0$ and $\Delta_t < 0$,  we obtain that  $\Delta_t $ increases by one with rate at least $\lambda + \lambda_0$ and decreases by one with rate at most $3 + 2\lambda$, and (\ref{liminf2}) easily follows.


\medskip

\noindent {\bf Remark.} In what follows we want to consider the unperturbed system as a random environment on which the
virtual particles will evolve.

\medskip

\noindent {\bf Step 2.}
From Proposition~\ref{prp:lln} and (\ref{liminf2})
it follows that  if $\lambda_0 > 3 + \lambda$, then
$$
\liminf_{t\to \infty} \frac1t(\tilde h_t(e_0) - \max\{\tilde h_t(e_{-1}) \tilde h_t(e_{1})\})> 0
$$
and, in particular there
exists $T_* < + \infty$ a.s.\ such that the number $\tilde n_0(s)$ of left virtual particles at the vertex $0$
and the number $\tilde n_1(s)$ of right virtual particles at vertex $1$ is positive, for all  $s > T_*$.
This also means that no virtual left (resp.\ right) particle created after $T_*$
enters $[2, + \infty)$ (resp.\  $(- \infty, - 1]$).

\medskip

\noindent {\bf Step 3.} If $\lambda_0 > 3 + \lambda$, then the total
number $ \sum_{k=2}^{+\infty} \tilde n_k(t)$ of virtual (in particular
virtual right) particles in $[2, + \infty)$ grows sub-linearly in $t$.
Indeed, if $\lambda_0 > 3 + \lambda$, then from Step 2 it follows that
there exists $T_* < + \infty$ a.s.\ such that $\Delta_s > 0$ for all $s
\ge T_*$. After time $T_*$  the system is essentially divided into two
half space systems. A virtual right-particle born after time $T_*$
(necessarily at site $1$) will then follow the path of a simple random
walk, possibly slowed down by the environment, until it first reaches
site $0$, when it will be annihilated; in particular, hitting
probabilities for virtual particles are given by a standard gambler's
ruin estimate.

Let $K$ be a positive integer. Consider all virtual
particles in $[2, + \infty)$ at large times $t > T_*$. They can be
subdivided into three classes:
\begin{itemize}
  \item  Those which were already present at time $T_*$, of which there
    are finitely many;
  \item  Those which were created after time $T_*$ and entered the
    half-line   $[K, +\infty)$ by time $t$. The number of such particles
    is almost surely asymptotically in $t$ bounded above by
    $(1+\varepsilon)t/K$, since the environment can possibly slow
    virtual particles down, but does not affect their trajectories.
  \item  Those which were created after time $T_*$, but never exited the
    box $[1, K]$. From Lemma 1 we get that the number of particles at
    any given site (except for $0$ and $1$) behaves sub-linearly ---
    notice that, since the proof of Lemma 1 is purely local, it still
    applies here to any site that is not directly affected by the
    additional nucleation at $e_0$, which is the case here. Thus,
    the number of such particles at any $x$ in the box is a.s.\ less
    than $t/K^2$ for all $t$ large enough.  In particular, the total
    number of such particles is smaller than $t/K$ for all $t$ large
    enough.
\end{itemize}
Thus, taking $K$ arbitrarily large we obtain the statement.

\noindent {\bf Step 4.} Observe, that for any $t < + \infty$ there
exists $0 < R(t) < +\infty$, random, such that $\tilde h_t (<x, x+1>)
=0$ if $|x| \ge R$. Thus, if $\lambda_0 > 3 + \lambda$, Step 3 implies
that $\tilde h_t (e)$ grows sub-linearly in time for all $e \neq e_0$.
Putting this together with Proposition~\ref{prp:lln} we get that $$
\lim_{t\to \infty} \frac{h'_t(e)}{t} =\lambda, \quad e \neq e_0, $$ for
such values of $\lambda_0$. By the monotonicity of the height dynamics,
this immediately extends to all $\lambda_0>0$.

\medskip

\noindent {\bf Step 5.} Applying the same line of reasoning as in Step 1 we see
that $\max\{\Delta_t,0\}$ grows at most sub-linearly if $\lambda_0 \leq 1$. Putting
together with Step 4, it follows that

$$
\lim_{t\to \infty} \frac{h'_t(e_0)}{t} =\lambda,  \quad \text{if } \lambda_0 \leq 1.
$$


\medskip

\noindent {\bf Step 6.} The
height process $h'_t(e_0)$ grows by one with rate at least
$\lambda+\lambda_0$,
and can only decrease if a particle (of the appropriate type) jumps
across the edge $e_0$ --- which happens with rate at most $1$. This implies that, if
$\lambda_0>1$, it is bounded below by a random walk with drift
$\lambda+\lambda_0-1$. On the other hand, the conclusion of Step 4
implies that for all $e \neq e_0$, $h'_t(e)$ satisfies a law of large
numbers with speed $\lambda<\lambda+\lambda_0-1$, thus for all $t$
large enough $n'_t (0)$ and $n'_t (1)$ are strictly positive, and so
$\Delta_t$ as well. Putting these two facts together and using the same
reasoning as above, we see that the increase and decrease rates for $h'_t (e_0)$
become exactly  $\lambda + \lambda_0$ and 1, respectively, for large $t$, which
implies the stated law for large numbers for $h'_t (e_0)$, if $\lambda_0>1$.

\medskip

\noindent This ends the proof of the main Theorem.

\bigskip

\noindent {\bf Acknowledgement.} The authors thank the
anonymous referees for many suggestions and for pointing out an error in the
previous version. We also thank prof. M. Myllys for the picture shown as Figure 1.
V.B. was partially supported by ANR Project Blan06-3-134462,
V.S. and M.E.V. were partially supported by Faperj grant E-26/100.626/2007 APQ1,
and CNPq Projeto Universal 48071/2006-1 and individual CNPq grants 300601/2005-0 and 302796/2002-9 respectively.

\bibliographystyle{spmpsci}
\bibliography{Biblio}

\end{document}